

\baselineskip=14pt
\parskip=10pt

\font\eighttt=cmtt8
\magnification=\magstephalf

\def\1{{\overline{1}}}
\def\2{{\overline{2}}}
\parindent=0pt
\overfullrule=0in

\def\frac#1#2{{#1 \over #2}}
\centerline
{\bf What is Mathematics and What Should it Be? }
\bigskip
\centerline
{\it Doron ZEILBERGER}

\qquad \qquad {\it Dedicated to Reuben Hersh on his 90th birthday}

{\bf Preamble}

In the classic ``{\bf The Mathematical Experience}'', Reuben Hersh and Philip Davis  gave a fresh look at mathematics and mathematicians, 
and showed us that notwithstanding Plato and Hilbert, mathematics is a {\it human} activity and culture,
and its pretensions to absolute truth are unfounded. Mathematics  {\it is what  mathematicians do}!.
This {\it leitmotif} was further expanded in the wonderful `sequel'
``{\bf Descartes' Dream: The World According to Mathematics}''  , also with Davis, and
the more philosophical {\it ``What Is Mathematics, Really?''}, written by Reuben all by himself.
The fact that mathematics is what (human) mathematicians do is further expounded,
beautifully, in ``{\bf Loving and Hating Mathematics}'', written in collaboration with Vera John-Steiner,
that, {\it inter-alia}, debunks G.H. Hardy's stupid quip that `Mathematics is a young men's game', by
describing many excellent women mathematicians, and many excellent `old' mathematicians.

What Hersh and Davis preached back in 1980 is still true today, but it has
to be tweaked quite a bit, since very soon (in $\leq 50$ years), mathematics would be
`what (machine) mathematicians do', and machines have different `emotional life'
than humans. Hence mathematics would finally, {\it hopefully}, become a true science.

{\bf What is Mathematics (as practiced today)?}

$\bullet$ A religion with its doctrines and dogmas

$\bullet$ A game with its (often arbitrary) rules

$\bullet$ An (intellectual) athletic competitive sport, akin to Chess and Go

$\bullet$ An art form with rigid rules

One thing it is not is a {\bf science}. Scientists, by definition, are trying to discover the {\it truth} about
the outside world. Mathematicians do not care about discovering the truth about the mathematical world.
All they care about  is playing their artificial game, called [rigorous] {\it proving}, and observing their strict dogmas.

{\bf  Beware of Greeks Carrying Mathematical Gifts}

Once upon a time, a long time ago, mathematics was indeed a true science, and its  practitioners devised
methods for solving practical mathematical problems. The ancient Chinese, Indian, and Babylonian mathematicians
were dedicated good scientists. Then came a major set-back, the Greeks!

The dirty (open) secret of the ``{\it enlightened}'' Athenian `democracy' was that it was a {\bf slave} society.
Since slaves did all the manual labor, the rich folks had plenty of time to
contemplate their navels, and to ponder about {\it the meaning of life}. Hence {\it Western Philosophy},
with its many {\it pseudo-questions} was developed in the hands of Plato, Aristotle, and their buddies,
and `Modern' pure mathematics was inaugurated in the hands of the gang of Euclid et. al.

Not all Greek mathematics was like that, and Archimedes, Heron, and many others also did excellent applied mathematics, but
even Archimedes was bound by the Euclidean `party line'.

Let's digress and summarize the pernicious Greek influence.

{\bf A Brief History of Mathematics as a Sequence of (Unsuccessfully!) Trying to Answer Stupid Questions}

$\bullet$ {\bf Stupid Question $1$}: Prove the Parallel Postulate, i.e. make it a theorem.

Many, very smart, people tried, in vain, to solve this problem,  until it turned out, as we all know, to be impossible.

This was really an {\it artificial puzzle}. Play the game of `logical deduction', starting with the (other) axioms,
and step-by-step construct an artificial structure called `rigorous mathematical proof', whose {\it bottom line} is
the statement that given a (geometrical) line and a point outside it, there is exactly one line through the point that does not meet the line.

$\bullet$  {\bf Stupid Question $2$}: Trisect an arbitrary angle only using straight-edge and compass.

Many, very smart people, tried in vain, to solve this problem,  until it turned out, in the 19th century,
to be {\bf impossible} 
(trisecting an angle involves constructing a cubic-algebraic number, while it is not too hard to show
that a {\bf necessary} condition for a ruler-and-compass constructible number is that its minimal equation is a power of $2$).

$\bullet$  {\bf Stupid Question $3$}: Double the cube. 

Ditto, $2^{1/3}$ is a cubic-algebraic number.

$\bullet$  {\bf Stupid Question $4$}: Square the Circle.

Many, very smart people, tried (and some still do!) 
in vain, to solve this problem,  until it turned out, in the late 19th century, thanks to Lindemann
(inspired by Hermite) to be {\bf impossible}.
Indeed, $\pi$, (and hence $\sqrt{\pi}$) is {\it transcendental},
all the more reason why it is not constructible.

Note that these three classical Delfian problems were really artificial puzzles, in an artificial game,
whose artificial rules are `only use straight-edge and compass'.

$\bullet$ {\bf Stupid Question $5$}: Find a general, closed-form, formula, as an expression in $p,q,r,s$,
{\it only using addition, multiplication, division, and root-extraction}
for a solution of the quintic equation
$$
x^5 +p x^3 + qx^2 +r x+ s=0 \quad .
$$

Recall that the analogous question  for a quadratic equation was answered, essentially, by the Babylonians.
The questions for the cubic and quartic  were answered by Tartaglia and Ferrari (and published by Cardano, ca. 1530), but
after that, many smart people, for about three hundred years, tried, {\it in vain}, to `solve' the quintic equation.

Let's digress to examine how the Renaissance mathematicians `solved' the cubic equation
$$
x^3+px+q=0 \quad .
$$
(Recall that we can always transform $x^3+a_2 x^2+ a_1 x+ a_0=0$ to the above form by writing $x=y-\frac{a_2}{3}$.)

We first do an {\it ad-hoc  trick}, writing  $x=u+v$, where $u$ and $v$ are {\bf to be determined}.

Then
$$
(u+v)^3+p(u+v)+q=0 \quad .
$$
Expanding
$$
u^3+3u^2v+3uv^2+v^3+p(u+v)+q=0 \quad .
$$
Rearranging:
$$
u^3+v^3+3u^2v+3uv^2+ p(u+v)+q=0 \quad .
$$
Replacing $3u^2v+3uv^2$ by $3uv(u+v)$ we get
$$
u^3+v^3+3uv(u+v)+ p(u+v)+q=0 \quad .
$$
Factoring out $(u+v)$ from the third and fourth terms:
$$
u^3+v^3+(3uv+p)(u+v)+q=0 \quad .
$$
We now do {\it wishful thinking},  {\bf demanding} that
$$
3uv+p=0 \quad .
$$
In other words
$$
uv= -\frac{p}{3} \quad .
$$
Going back to the above equation we have
$$
u^3+v^3+0+q=0 \quad .
$$
So
$$
u^3+v^3= -q \quad .
$$
Cubing the equation  $uv=-\frac{p}{3}$, we get
$$
u^3v^3= -\frac{p^3}{27} \quad .
$$
Hence the sum of $u^3$ and $v^3$ is $-q$ and their product is
$ -\frac{p^3}{27}$, hence {\bf both} are solutions of the
{\bf quadratic equation}
$$
z^2+qz -\frac{p^3}{27}=0 \quad .
$$

By the quadratic formula, the two roots $u^3,v^3$ are
$$
\frac{-q \pm \sqrt{q^2+ \frac{4}{27}p^3}}{2}
$$
So one root of the cubic is
$$
\left (\frac{-q + \sqrt{q^2+ \frac{4}{27}p^3}}{2} \right )^{\frac{1}{3}}+
\left (\frac{-q - \sqrt{q^2+ \frac{4}{27}p^3}}{2} \right )^{\frac{1}{3}} \quad .
$$

The other two roots are $\omega u + \omega^2 v$ and $\omega^2 u + \omega v$, where $\omega$ is a root of $\omega^2-\omega+1=0$.

So by an {\it ad hoc} trick, and {\bf a lot of luck}, we {\it reduced} solving a cubic to that of solving
a quadratic.

Similar ad-hoc tricks work for the quartic, but mathematicians were stumped for three hundred years,
until, famously, Ruffini, Abel, and most notably, Galois, proved that it is {\bf impossible}.

In hindsight, `solvability by radicals' is an {\it artificial game} with an {\it artificial} set of `legal moves', and
Galois et. al. proved that it is impossible to get from position $A$, the quintic equation, to position $B$,
an algebraic formula only using addition, multiplication, division, and root-extraction.

I admit that while the question itself `solve a quintic by radicals', turned out (at least in hindsight) to be
very stupid, it lead to something not quite as stupid, {\it Group theory}, and {\it Galois theory}.

$\bullet$  {\bf Stupid Question $6$}: Find a  ``rigorous'' foundation to the non-rigorous Differential and Integral Calculus
of Newton and Leibnitz, thereby `resolving' the paradoxes of Zeno, and `addressing' Bishop Berkeley's critique.

This was {\it allegedly} `solved' by Cauchy and  Weierstrass, but their `solution' was unnecessarily complicated
and pedantic, creating so-called `Real Analysis', one of the most unattractive courses in the
math-major curriculum, where one does scholastic mental gymnastics to `prove' intuitively obvious facts.

As I pointed out in my essay
{\it `` `Real' Analysis is a Degenerate Case of Discrete Analysis''},

[{\eighttt ``New Progress in Difference Equations'', edited by Bernd Aulbach, Saber Elaydi, and Gerry Ladas, (Proc. ICDEA 2001), Taylor and Francis, London}]
available on-line from my website 

{\tt http://www.math.rutgers.edu/\~{}zeilberg/mamarim/mamarimPDF/real.pdf} \quad ,

a much better, {\bf conceptually simpler}, and more honest, foundation, is to {\bf chuck} infinity, and {\it limits}
altogether, and replace the derivative of a function $f(x)$ that is defined as the limit of $(f(x+h)-f(x))/h$ as $h$ goes to zero,
by the difference operator $\Delta_h f(x):=(f(x+h)-f(x))/h$, where $h$ is a very small, but {\bf not} `infinitesimal'.
Since $h$ is so tiny (and unknowable) it is more convenient to  leave it {\it symbolic}
(analogously to the way physicists write $h$ for Plank's constant).
It is true that technically things get a bit messier, for example, the product rule is a little more complicated, but
this is a small price to pay, and at the end of the day, when you replace the symbolic $h$ by
$0$ you get the familiar rules.

We live in a {\bf finite} and {\bf discrete} world, and the infinite and the continuous are mere {\it optical illusions}.

Ironically, the way that `continuous' differential equations are (numerically) solved today, by computers, is by
{\it approximating} them with finite-difference equations, and numerical analysts make a living by proving
{\it a priori} {\it error-estimates}. The true equations are finite-difference equations to begin with, but
the mesh-size is too small (and unknowable, the above mentioned $h$), that it would be impractical to solve
numerically, and one has to replace it by a far coarser grid.

$\bullet$  {\bf Stupid Question $7$}: Is there a set whose cardinality is strictly between that of the integers and that
of the real numbers?

This is Hilbert's first problem, the so called {\it Continuum Hypothesis}, or CH for short.
It is really a {\it pseudo-question}, since it pertains to two `infinite', and hence fictional, sets.
It is really an artificial puzzle. Can you reach, by finitely many legal moves (in the game called `logical deduction')
starting with the axioms of ZFC, either the statement that there exists such a set, or the statement that
there does not exist such a set. Paul Cohen famously proved that neither! 
The conventional way of saying it is that CH
is {\bf independent} of ZFC, but what Cohen {\it really} (brilliantly!) meta-proved was that the question was {\it stupid}!
Neither CH, nor its negation, are reachable in the logical deduction game starting from the axioms of ZFC, hence
{\it both} CH, and {\it ZFC}, are devoid of content, and what Cohen (and G\"odel) meta-proved was
that the so-called {\it infinity} does not exist! [taken at {\it face value}].

$\bullet$  {\bf Stupid Question $8$}: Are the axioms of arithmetic consistent?

This, Hilbert's second problem, was famously shown to be a stupid question by Kurt G\"odel.
The conventional, polite, interpretation is that there exist {\bf undecidable} statements,
but a more honest interpretation is that every statement that involves quantifies over
`infinite' sets is {\it a priori meaningless}. Some such statements, e.g. that $n+1=1+n$ for
{\it every} natural number, can be made {\it a posteriori meaningful}, by thinking of $n$ 
as a {\it symbol}, but what G\"odel proved was that many such statements may not be
resurrected like that.

$\bullet$  {\bf Stupid Question $9$}: ``Given a diophantine equation with any number of unknown quantities
and with rational integral numerical coefficients: {\it To devise a process according to which it can be determined
by a finite number of operations whether the question is solvable in rational integers}.''

This was famously (brilliantly) meta-proved to be a stupid question by Yuri Matiyasevich, standing on the shoulders
of Martin Davis, Hillary Putnam, and Julia Robinson. What he really proved was that, while for
some {\it specific} diophantine equations, it is possible to prove that there are {\bf no solution}, (e.g.
$x^2-2y^2=0$ and the equation $x^n+y^n=z^n \, , \, n>2 \, , \, xyz \neq 0$), 
while for other ones it is possible to find  many (in fact, {\it potentially} infinitely many) solutions
(e.g. $x^2-2y^2=1$), many such 
(seemingly sensible, but in fact meaningless) questions can't be resolved in the artificial game called
proving number-theoretic statements. In other words, there are not only {\it a priori} meaningless
(as all statements involving `infinite' sets are), they are also {\it a posteriori} meaningless.

{\bf Many Games Have Unreachable Positions}

Time and time again, mathematicians realized that, in their human-made, artificial game, that they naively
believed to be the {\it real thing}, many things are {\bf impossible}. Sometimes it is easy to prove the
{\it impossibility}, for example tiling, by domino-pieces, an $n \times n$ checker-board where two opposite corners
have been removed, when $n$ is odd.
Since $n^2-2$ is odd and any covering by domino-pieces
must cover an even number of unit-squares. Other times the impossibility  is more subtle, this
is the case When $n$ is even
(e.g. $n=8$, the usual chess board).
Coloring the unit-squares (as in a chess board)
white and Black alternatively, any putative domino-covering must cover an equal number of
black and white unit-squares, but if you remove two opposite corners (from an $n \times n$ board with $n$ even),
they must be of the same color, hence the remainder has an excess of two of the other color.

A less obvious example of proven impossibility is in the solitaire game called the {\it fifteen puzzle}, that Sam Loyd
was safe in offering a large prize for its solution (transposing the $14$ and $15$ and leaving the rest the same),
since, as proved in one of the first volumes of  James Joseph Sylvester's periodical, {\it American Journal of Mathematics},
is {\bf impossible}.

{\bf Today's Mathematics Is a Religion}

Its central dogma is {\bf thou should prove everything rigorously}. Let me describe two examples of the religious
fanaticism of two of my good friends, George Andrews and Christian Krattethaler. I admire their mathematics,
and I like them as people, but I was disappointed (but also impressed) at their {\it fanaticism} 
(and [misguided] {\it integrity}).

Drew Sillls and I {\bf disproved} (in the everyday sense of the word), {\it conclusively}, a long-standing
conjecture by Hans Rademacher (one of the greatest number theorists of the 20th century, and
George Andrews' thesis advisor). Since the disproof was not `rigorous', George Andrews refused to consider it
for the {\it Proc. of the National Acad. of Sci.} (PNAS), and when we submitted it, it was rejected, because he
was consulted. Then we went to the less prestigious journal {\it Mathematics of Computation},
and once again it was rejected, since it was `too computational'. Then we went `down'  to 
the Journal ``Experimental Mathematics'', and once again it was rejected, by its editor,
Yuri Tschinkel, since it was `too experimental'. Luckily (for Drew Sills, I already was fully
promoted, and would have left it in arxiv.org) it was finally accepted in another journal.

Another amusing example of the religious dogmatism of current, otherwise reasonable, mathematicians,
is when my collaborators, Manuel Kauers, Christian Koutschan, and I submitted to the
not-exactly-prestigious journal {\it S\'eminaire Lotharingien de Combinatoire} an article

{\it A Proof of George Andrews' and Dave Robbins' q-TSPP Conjecture (modulo a finite amount of routine calculations)}

It contained a proof-plan, that, given enough computer resources would give a fully rigorous proof,
and we had plenty of evidence that the plan should work. But Christian Krattenthaler objected to the title!
He would have been happy to accept it if we changed the title to:

``{\it A Proposal for a Possible Computer Proof of George Andrews' and Dave Robbins' q-TSPP Conjecture}'',

but we refused, so it remained in our websites (and the arxiv). Luckily, a year later, we found
a way to prove it with today's computers, and this, {\it fully rigorous} version was
gladly accepted by George Andrews for the PNAS. But all the main ideas were already contained in the
previous, semi-rigorous version, with clear evidence that given enough computing power, it should work.

{\bf Today's Mathematics is a Competitive Sport}

Number Theory is full of conjectures that are obviously true, but humans are unable to prove them
[in their outdated, narrow-minded sense of the word, meaning `rigorous-proof'].
For example, Yitang Zhang got instantly famous by getting ever-so-close to the twin-prime conjecture.
But even the still-open twin-prime conjecture is so far from what is definitely true, by
very plausible heuristics. The twin-prime conjecture asserts that there are infinitely primes $p$ such
that also $p+2$ is prime. The true fact is that there are lots and lots of them, not just `infinitely many'!
In fact, there are $O(n/\log(n)^2)$ such twin-prime pairs less than $n$. But mathematicians do not care about
{\it truth}, they only care about playing their (artificial!) game.

{\bf Today's Mathematics is an Artificial Game}

As I demonstrated at length above.

{\bf Today's Mathematics is an Art Form}

Paul Erd\H{o}s famously talked about {\it proofs from the book}, and, indeed, {\it elegance}, 
and {\it beauty} are  great values cherished by human mathematicians. 

Humans, starting with Plato, via Dirac, Hardy, and almost everyone else,
waxed eloquently about {\it beauty}. G. H. Hardy went as far as saying: 

{\it ``There is no place for ugly mathematics''}.

Excuse me, my dear Hardy, this is even stupider than your `young man's game' unfortunate quip [debunked by Hersh; btw, note the word {\it game}].
Beauty is only skin-deep,  and is also in the {\it eyes of the beholder}, and who are you to {\it exclude}
`ugly' mathematics! Ironically, personally, I find most of your hard-analysis
mathematics much uglier than the average mathematics, but this is besides the point.

In a cute New York Times  [April 16, 2017] Opinion piece,
entitled `{\it Beautiful Equations}', math-groupie and distinguished Cornell psychiatry professor, Richard A. Friedman
describes a psychological experiment conducted
by Semir Zeki (in collaboration with three other co-authors: John Paul Romaya, Dionigi M. T. Benincasa, and guess who?,
Abel-laureate and Fields-medalist  Sir Michael Atiyah), 

[{\it `The experience of mathematical beauty and its neural correlates'},
Front. Hum. Neurosci., 13 February 2014 {\tt https://doi.org/10.3389/fnhum.2014.00068}] \quad .

In that experiment, fifteen mathematicians were asked to rank a list of sixty famous equations according
to beauty, and then electrodes were connected to their brains (using fMRI scanners), 
and there turned out to be a high correlation between the subject's beauty-rank of the formula before
the scan, and the activation of the area of the brain responsible for aesthetic pleasure.

This is a very interesting psychology experiment, that tells a lot about humans, but does not say
anything about what mathematics {\it should} be.

It is instructive to see what rankded first and what ranked last. The number-one-hit-formula was Euler's
$$
e^{i\pi}+1 \, = 0 \quad,
$$
and the bottom-one, the {\it ugliest} (in the eyes of the mathematicians guinea-pigs) was
Ramanujan's formula for $\frac{1}{\pi}$
$$
{{1} \over {\pi}} \, = \,
2 \sqrt{2} \, \sum_{k=0}^{\infty} \,
\frac{(1103+26390k) \,(1/4)_k (1/2)_k (3/4)_k}{k!^3} \cdot \frac{1}{ 99^{4k+2}} \quad ,
$$
where $(a)_k=a(a+1) \dots (a+k-1)$.

This shows that human mathematicians are superficial. In the `eyes of God' Ramanujan's formula is
much prettier than Euler's. It is a deep relationship enabling a very fast computation of
$\pi$ to billions of decimals. In contrast, Euler's {\it beauty-queen} is an utter triviality,
only a bit less trivial than
$$
1+1=2 \quad .
$$
In fact, these two equations have something in common, they may be viewed as {\it definitions}.
$2$, by definition, is $1+1$, and $\pi$, by definition, is the smallest real number larger than $0$
for which $e^{ix}$ happens to be equal to $-1$. In other words, the smallest $x$
for which $\sin x$ equals $0$ (and hence $\cos x$ equals $-1$, and hence $e^{ix}=\cos x+ i \sin x= - 1+ i \cdot 0= -1$).

While there is no harm in enjoying  mathematics for its (subjective!) beauty, it is wrong
to make it a {\it defining} property, and to {\it exclude} what one (subjectively!) finds {\it ugly}.
Every statement that has a {\it proof from the book}, is {\it ipso facto}, {\it trivial}, at least,
{\it a posteriori}, since all deep statements have long proofs. In my eyes the most beautiful theorems
are those with succinct statements for which the shortest known (and hopefully any) proof
is very long. So in my eyes, both the Appel-Haken {\it Four Color Theorem}, and Tom Hales' {\it Kepler's ex-Conjecture}
should be in the {\bf book}.

{\bf Today's Mathematics is an Elitist and Exclusive Club}

It has a very intricate {\it social structure}, with a fairly rigid pecking order, and
hierarchy, and `peer-reviewed journals',
with pompous, often {\it sadistic}, editors, who must enjoy rejecting submissions
(or else they would refuse to take the job).

Some areas are considered very {\it respectable}, while other ones are {\it slums}.
This is a dynamic process, and areas come and go out of fashion, but there are always,
like in a high school, the `cool kids lunch table', and there are always the outcasts and
the {\it pariahs}. Humans will be humans!

{\bf Today's Mathematics is Not a Science}

As we saw above, mathematics today is many things, but one thing that it is {\bf not}, is a {\it science}.
It is amazing, that, nevertheless, mathematics was so {\it effective} in science, so convincingly
told by Eugene Wigner.  One reason it was so effective was that the kind of mathematics that
scientists needed was either discovered or rediscovered by themselves (e.g. Heisenberg
rediscovered matrix algebra, {\it ab initio}), and they develop their own brand of mathematics,
{\bf without} the mathematicians' misplaced obsession with {\it rigor} (e.g. Quantum Field Theory and
the Renormalization Group), and the great success of mathematics in science is
{\bf in spite} of the mathematicians' superstitious dogmas. Imagine how more
effective it would have been if they threw away their artificial {\it shackles}.

{\bf The Computer Revolution}

The reason mathematics today is the way it is, is due to the {\it contingent} fact that it was
developed without computers. Most of the questions (and pseudo-questions) that occupied
mathematicians through the ages are now {\it moot} and {\it irrelevant}. It may
be a good idea to start {\it all over}, and develop mathematics {\it ab initio},
not peeking at the human-generated mathematics done so far, and taking full advantage
of computers.

{\bf Let's Make Mathematics a Science Again!}

Mathematics {\bf should} become a {\it science}, and its main {\it raison-d'\^etre} should be
the  discovery of mathematical {\bf truth (broadly defined!)}. In particular, one
should abandon the {\it dichotomy} between {\it conjecture} and {\it theorem}.

If a mathematical question is in doubt (and it looks that it can go either way: true, false (or meaningless (i.e. `undecidable')),
then it is a mathematical {\it question}. If there is overwhelming empirical and/or heuristic evidence, then
it is a {\it theorem}. What was formerly called a theorem should be renamed 
`rigorously-proved statement' [using the artificial game of logical deduction].

{\bf Final Words}

Don't get me wrong, while you may call me a {\it self-hating-mathematician}, I do {\it love} elegant proofs,
and still enjoy traditional standards. So like  Reuben Hersh's wonderful book mentioned above, I
{\bf Love and Hate} (traditional) mathematics. But I {\bf do} believe that it is time to make it
a {\bf true} science.

\bigskip
\hrule
\bigskip
Doron Zeilberger, Department of Mathematics, Rutgers University (New Brunswick), Hill Center-Busch Campus, 110 Frelinghuysen
Rd., Piscataway, NJ 08854-8019, USA. \hfill\break
Email: {\tt DoronZeil at gmail dot com}   \quad .

\end